\documentclass[11pt,a4paper,twoside]{article}
\usepackage{amssymb}
\usepackage{amsmath}
\usepackage{fancyhdr}
\usepackage{amsthm}
\usepackage[cp1250]{inputenc}

\title{The classification on simple Moufang loops}
\date{}
\author{SANDU N. I.}
\begin{document}
\maketitle

\begin{abstract}
Let $C(F)$ be a matrix Cayley-Dickson algebra over field $F$. By
$M_0(F)$ we denote the loop containing of all elements of algebra
$C(F)$ with norm $1$. It is shown in this paper that with
precision till isomorphism the loops $M_0(F)/<-1>$ they and only
they are simple non-associative Moufang loops, where $F$ are
subfields of algebraic closed field.
\smallskip\\


\end{abstract}
\bigskip

The purpose of this paper is to classify the non-associative
simple Moufang loops. The Moufang loop is simple if it has no
non-trivial proper normal subloops, or equivalently, if it has no
non-trivial proper homomorphic images. For basic definitions and
properties of Moufang loops see [1].

It is well known that for an alternative algebra $A$ with the unit
element $1$  the set $U(A)$ of all invertible elements of $A$
forms a Moufang loop with respect to multiplication [2]. Let now
$C(F)$ be a matrix Cayley-Dickson algebra over arbitrary field
$F$, and let $M_0(F)$ be the set of all elements of $C(F)$ with
norm $1$. $M_0(F)$ is a normal subloop of $U(C(F))$. Let
$Z(M_0(F))$ be the center of $M_0(F)$. Paige L. shows in [3,
Theorem 4.1 and corollary to Lemma 3.4] that $M_0(F)/Z(M_0(F))$ is
a simple, non-associative, Moufang loop and center $Z(M_0(F))$ is
a group of order $2$ if the characteristic of $F$ is not $2$;
otherwise $Z(M_0(F)) = 1$. Further, using this result and the
powerful apparatus of finite groups theories, in [4] Liebeck M.
finalizes the classification of finite non-associative simple
Moufang loops, started in [5, 6]. He shows that such loops are
isomorphic to one of the loops $M(GF(q))$, where $GF(q)$ is a
finite Galois field, modulo the center. Within this paper all
non-associative simple Moufang loops are classified through
methods of alternative algebras. With precision till isomorphism
it is one of loops $M_0(C(F))/<-1>$, where $C(F)$ denotes a matrix
Cayley-Dickson algebra $C(F)$ over a subfield $F$ of algebraically
closed field. Further, the needed results of alternative algebras
from [7, 8] will be used without reference.

By analogy to Lemma 1 from [9] it is proved.
\smallskip\\

\textbf{Lemma 1.} \textit{Let $A$ be an alternative algebra and
let $Q$ be a subloop of $U(A)$. Then the restriction of any
homomorphism of algebra $A$  upon $Q$ will be a homomorphism on
the loop. More concretely, any ideal $J$ of $A$ induces a  normal
subloop $Q \cap (1 + J)$ of $Q$.}
\smallskip\\

Let $L$ be a free Moufang loop, let $F$ be a field and let $FL$ be
a loop algebra of loop $L$ over field $F$. We remind that $FL$ is
a free module with basis $\{g \vert g \in L\}$ and the
multiplication of elements of the basis is defined by their
multiplication in loop $L$. Let $(u,v,w) = uv\cdot w - u\cdot vw$
denote the associator of elements $u, v, w$ of algebra $FL$. We
denote by $I$ the ideal of loop algebra $FL$, generated by the set

$$\{(a,b,c) + (b,a,c), (a,b,c) + (a,c,b) \vert \forall a, b, c \in L\}.$$ It
is shown in [9] that algebra $FL/I$ is alternative and loop $L$ is
embedded (isomor\-phically) in the loop $U(FL/I)$. Further we
identify the loop $L$ with its isomorphic image in $U(FL/I)$.
Hence the free loop $L$ is a subloop of loop $U(FL/I)$. Without
causing any misunderstandings, like in [9], we will denote by $FL$
the quotient algebra $FL/I$ and call it $''$loop algebra$''$ (in
inverted commas). Sums $\sum_{g\in L}\alpha_gg$, are elements of
algebra $FL$, where $\alpha_g \in F$. Further, we will identify
the field $F$ with subalgebra $F1$ of algebra $FL$, where $1$ is
the unit of loop $L$.

Let now $Q$ be an arbitrary Moufang loop. Then $Q$ has a
representation as a quotient loop $L/H$ of the free Moufang loop
$L$ by the normal subloop $H$. We denote by $\omega H$ the ideal
of $''$loop algebra$''$ $FL$, generated by the elements $1 - h$
($h \in H$). By Lemma 1 $\omega H$ induces a normal subloop $K = L
\cap (1 + \omega H)$ of loop $L$ and $F(L/K) = FL/\omega H$.

We denote $L/K = \overline{Q}$, thus $FL/\omega H =
F\overline{Q}$. As every element in $FL$ is a finite sum  $\sum_{g
\in L}\alpha_gg$, where $\alpha_g \in F$, $g \in L$, then the
finite sum $\sum_{q \in \overline{Q}}\alpha_qq$, where $\alpha_q
\in F$, $q \in \overline{Q}$ will be elements of algebra
$F\overline Q$.  Let us determine the homomorphism of $F$-algebras
$\varphi: FL \rightarrow F(L/H)$ by the rule
$\varphi(\sum\lambda_qq) = \sum\lambda_qHq$. The mapping $\varphi$
is $F$-linear, then for $h \in H$, $q \in L$ we have $\varphi((1 -
h)q) = Hq - H(hq) = Hq - Hq = 0$. Hence $\omega H \subseteq
\text{ker}\varphi$. The loop $\overline{Q}$ is a subloop of loop
$U(F\overline{Q})$ and as $\omega H \subseteq \text{ker}\varphi$,
then the homomorphisms $FL \rightarrow FL/\omega H = F\overline Q$
and $FL \rightarrow FL/\text{ker}\varphi = F(L/H) = FQ$ induces a
homomorphism $\pi$ of  loop $\overline{Q}$ upon loop $Q$. Hence we
have.
\smallskip\\

\textbf{Lemma 2.} \textit{Let $Q$ be an arbitrary Moufang loop.
Then the loop $\overline{Q}$ is embedded in loop of invertible
elements $U(F\overline{Q})$ of alternative algebra $F\overline{Q}$
and the homomorphism $L \rightarrow FL/\omega H$ of $''$loop
algebra$''$ $FL$ induces a homomorphism $\pi: \overline Q
\rightarrow Q$ of loops.}
\smallskip\\

Let now $Q$ be a simple Moufang loop. Then $\text{ker}\pi$ will be
a proper maximal normal subloop of $\overline Q$. Let $J_1,$ $J_2$
be  proper ideals of algebra $F\overline Q$. We prove that the sum
$J_1 + J_2$ is also a proper ideal of $F\overline Q$. Indeed, by
Lemma 1 $K_1 = \overline Q \cap (1 + J_1),$ $K_2 = \overline Q
\cap (1 + J_2)$ will be  normal subloops of loop $\overline Q$. We
have that $K_1 \subseteq \text{ker}\pi$, $K_2 \subseteq
\text{ker}\pi$. Then product $K_1K_2 \subseteq \text{ker}\pi$, as
well. But $K_1K_2 = (\overline Q \cap (1 + J_1))(\overline Q \cap
(1 + J_2)) = \overline Q \cap (1 + J_1)(1 + J_2) = \overline Q
\cap (1 + J_1 + J_2 + J_1J_2) = \overline Q \cap (1 + J_1 + J_2)$.
Hence $\overline Q \cap (1 + J_1 + J_2) \subseteq \text{ker}\pi$,
i.e. $\overline Q \cap (1 + J_1 + J_2)$ is a proper normal subloop
of $\overline Q$. Then from Lemma 1 it follows that $J_1 + J_2$ is
a proper ideal of algebra $F\overline Q$, as required.

We denote by $S$ the ideal of algebra $F\overline  Q$, generated
by all proper ideals $J_i$ ($i \in I$) of $F\overline Q$. Let us
show that $S$ is also a proper ideal of algebra $F\overline Q$. If
$I$ is a finite set, then the statement follows from first case.
Let us now consider the second possible case. The algebra
$F\overline Q$ is generated as a $F$-module by elements $x \in
\overline Q$. Let there be such ideals $J_1, \ldots, J_k$ that for
element $1 \neq a \in \overline Q$ $a \in \sum J_i$ and let us
suppose that for element $b \in \overline Q$ $b \notin \sum J_i$.
We denote by $T$ the set of all ideals of algebra $F\overline Q$,
containing the element $a$, but not containing the element $b$. By
Zorn's Lemma there is a maximal ideal $I_1$ in $T$. We denote by
$I_2$ the ideal of algebra $F\overline  Q$, generated by all
proper ideals of $F\overline  Q$ that don't belong to ideal $I_1$.
Then $S = I_1 + I_2$. $I_1, I_2$ are proper ideals of $F\overline
Q$ and by first case $S$ is also proper ideal of $F\overline Q$.
By Lemma 1 $K = \overline Q \cap (1 + S)$ is a normal subloop of
$\overline Q$. We denote $\overline{\overline Q} = \overline Q/K$.
Then $F\overline{ \overline Q} = F\overline Q/S$ is a simple
algebra. As $K \subseteq \text{ker}\pi$ then $\pi$ induce a
homomorphism $\rho: \overline{\overline Q} \rightarrow Q$. Hence
we prove.
\smallskip\\

\textbf{Lemma 3.} \textit{Let $Q$ be a simple non-associative
Moufang loop. Then $F\overline{\overline Q}$ is a simple
alternative algebra and the homomorphism $\pi: \overline Q
\rightarrow Q$ induces a homomorphism $\rho: \overline{\overline
Q} \rightarrow Q$.}
\smallskip\\

Let $F$ be an arbitrary field. Let us consider a classical matrix
Cayley-Dickson algebra $C(F)$. It consists of matrices of form
$\left(
\begin{array}{ll} \alpha_1 & \alpha_{12}\\
\alpha_{21} & \alpha_2  \end{array} \right)$, where $\alpha_1,
\alpha_2 \in F$, $\alpha_{12}, \alpha_{21} \in F^3$. The addition
and multiplication by scalar of elements of algebra $C(F)$ is
represented by ordinary addition and multiplication by scalar of
matrices, and the multiplication of elements of algebra $C(F)$ is
defined by the rule

$$ \left(
\begin{array}{ll} \alpha_1 & \alpha_{12}\\
\alpha_{21} & \alpha_2  \end{array} \right) \left(
\begin{array}{ll} \beta_1 & \beta_{12} \\ \beta_{21} & \beta_2
\end{array} \right) = $$
$$ \left( \begin{array}{ll} \alpha_1 \beta_1 + (\alpha_{12},
\beta_{21}) & \alpha_1 \beta_{12} + \beta_2 \alpha_{12} -
\alpha_{21} \times \beta_{21} \\ \beta_1 \alpha_{21} + \alpha_2
\beta_{21} + \alpha_{12} \times \beta_{12} & \alpha_2 \beta_2 +
(\alpha_{21}, \beta_{12}) \end{array} \right), $$ where for
vectors $\gamma = (\gamma_1, \gamma_2, \gamma_3),$ $\delta =
(\delta_1, \delta_2, \delta_3) \in A^3$  $(\gamma, \delta) =
\gamma_1\delta_1 + \gamma_2\delta_2 + \gamma_3\delta_3$ denotes
their scalar product and $\gamma \times \delta = (\gamma_2\delta_3
- \gamma_3\delta_2, \gamma_3\delta_1 - \gamma_1\delta_3,
\gamma_1\delta_2 - \gamma_2\delta_1)$ denotes the vector product.
Algebra  $C(F)$ is alternative. It is also  quadratic over $F$,
i.e. each element $a \in C(F)$ satisfies the identity

$$a^2 - t(a)a + n(a) = 0, n(a), t(a) \in F$$ and admits
composition, i.e.

$$n(ab) = n(a)n(b)$$ for $a, b \in C(F)$. Track $t(a)$ and norm
$n(a)$ are defined by the equalities $t(a) = \alpha_1 + \alpha_2,$
$n(a) = \alpha_1\alpha_2 - (\alpha_{12}, \alpha_{21})$.

We have $M_0(F) = \{u \in C(F) \vert n(u) = 1\}$. Further, $n(1) =
1$, and it follows from the relations $n(ab) = n(a)n(b)$,
$n(\alpha a) = \alpha^2n(a)$ that $-1 \in M_0(F)$. Obviously $-1$
belongs to the center of algebra $C(F)$. Then $-1$ belongs to the
center $Z(M_0(F)$ of loop $M_0(F)$. Therefore the subloop $<-1>$,
generated by element $-1$, is normal in $M_0(F)$ and from Paige's
results [3], presented at the beginning of the article it follows.
\smallskip\\

\textbf{Lemma 4.} \textit{Let $F$ be an arbitrary field. Then the
Moufang loop $M(F) = M_0(F)/<-1>$ of the matrix Cayley-Dickson
algebra $C(F)$ is simple and the loop $<-1>$ coincide with center
of loop $M_0(F)$}.
\smallskip\\

By Lemma 4 the center $Z$ of loop $M_0(F)$ coincides with subloop
$<-1>$. As $M_0(F)/Z$ is a simple loop, a question appears. Is the
center $Z$ of loop $M_o(F)$ emphasized by the direct factor? The
answer is negative. Let field $F$ consist of 5 elements and let
$H$ be a direct completion of center $Z$. If $\alpha$ were the
generator of the multiplicative group of field $F$, then one of
the elements $\pm\left(
\begin{array}{ll} \frac{1}{\alpha} & 0 \\ 0 & \alpha \end{array}
\right)$ would lie in $H$. The square of this element is equal to
$-1$, i.e., it lies in the intersection $H \cap Z$, which is
impossible. Therefore center $Z$ cannot have a direct factor in
$M_0(F)$.

Let now $P$ be an algebraically closed field and let $Q$ be a
simple non-associative Moufang loop. By Lemma 3 the loop
$\overline{\overline Q}$ is embedded in loop of invertible
elements of simple alternative algebra $P\overline{\overline Q}$.
We denote $\overline{\overline Q} = G$.

If $a \in G$, then it follows from the equality $aa^{-1} = 1$ that
$n(a)n(a)^{-1} = 1$, i.e. $n(a) \neq 0$. Associator $(a,b,c)$ of
elements $a, b, c$ of an arbitrary loop is defined by the equality
$ab\cdot c = (a\cdot bc)(a,b,c)$. Identity $(xy)^{-1} =
y^{-1}x^{-1}$ holds in Moufang loops. Therefore, if $a, b, c$ are
elements of Moufang loop $G$, then $u = (a,b,c) = (a\cdot
bc)^{-1})(ab\cdot c) = (c^{-1}b^{-1}\cdot a^{-1})(ab\cdot c)$,
$n(u) = n(c^{-1})n(b^{-1})n(a^{-1})n(a)\cdot n(b)n(c) =
n(c)^{-1}n(b)^{-1}n(a)^{-1}n(a)n(b)n(c) = 1$, i.e. $u \in M_0(P)$.
We denote by $G'$ the subloop generated by all associators of
Moufang loop $G$. If $G' = G$, then $G \subseteq M_0(P)$, i.e. the
loop $G$ is embedded in $M_0(P)$. Now we suppose that $G' \neq G$.
It is shown in [10, 11] that the subloop $G'$ is normal in $G$.
The finite sum $\sum_{g \in G}\alpha_gg$, where $\alpha_g \in P$,
$g \in G$ are elements of algebra $PG$.  Let $\eta :PG \rightarrow
P(G/G')$ be a homomorphism of $P$-algebras determined by rule
$\eta(\sum \alpha_gg) = \sum \alpha_ggG'$ ($g \in G$) and let
$P(G/G') = PG/\text{ker}\eta$. As the quotient loop $P(G/G')$ is
non-trivial, then $PG/\text{ker}\eta \neq PG$. Hence
$\text{ker}\eta$ is a proper ideal of $PG$. The algebra $PG$ is
simple. Then the ideal $\text{ker}\eta$ cannot be the proper ideal
of $PG$. Hence the case $G' \neq G$ is impossible and,
consequently, the loop $G$ is embedded in loop $M_0(P)$.

The alternative algebra $PG$ is simple. By Kleinfeld Theorem [12,
see also 7, 8] it is a Cayley-Dickson algebra over their center.
Field $P$ is algebraically closed. Then algebra $PG$ is split. The
matrix Cayley-Dickson algebra $C(P)$ is also split. But any two
split non-associative composition algebras over an algebraically
closed field are isomorphic. Therefore algebra $PG =
P\overline{\overline Q}$ is isomorphic to the matrix
Cayley-Dickson algebra $C(P)$.

If $M$ is a $F$-module, and $N$ is its subset, then the denotation
$F\{N\}$ means the $F$-submodule generated by $N$. We denote
$P\overline{\overline Q} = C_P(\overline{\overline Q})$.
Consequently, it is proved.
\smallskip\\

\textbf{Lemma 5.} \textit{Let $P$ be an algebraically closed field
and let $Q$ be a simple Moufang loop. Then loop
$\overline{\overline{Q}}$ is embedded in split Cayley-Dickson
algebra $C_P(\overline{\overline{Q}}) =
P\{\overline{\overline{Q}}\}$.}
\smallskip\\

Let $P$ be an algebraically closed field and let $H$ be a subloop
of loop $M_0(P)$. Let $a_{ij} = (a_{ij}^{(1)}, a_{ij}^{(2)},
a_{ij}^{(3)})$.   If the matrices elements $\alpha_i,
\alpha_{ij}^{(k)}$ of all
matrices $\left(\begin{array}{ll} \alpha_1 & \alpha_{12} \\
\alpha_{21} & \alpha_2 \end{array} \right)$ $ \in H$ generate the
subfield $F$ of field $P$, then we will say that $H$ is a
\textit{loop over  field} $F$. If $H$ is strictly contained into
$M_0(F)$, then we say that $H$ is a \textit{proper over field $F$
subloop} of loop $M_0(F)$.
\smallskip\\

\textbf{Lemma 6.} \textit{Let $F$ be an arbitrary field. Then the
Moufang loop $M_0(F)$ doesn't contain proper over $F$
non-associative subloops of type $\overline{\overline{Q}}$,
considered in Lemma 5.}
\smallskip\\

\textbf{Proof.} Let $P$ be an algebraic closing of field $F$. Let
us suppose the contrary, that loop $M_0(F) = L$ contains a proper
over $F$ non-associative subloop $H \subset L$ of type
$\overline{\overline{Q}}$. Let $P\{H\}$, $P\{L\}$ are  the matrix
Cayley-Dickson algebras, considered in Lemma 5. We consider the
subalgebras $C_F(H) = F\{H\} \subseteq P\{H\}$, $C_F(L) = F\{L\}
\subseteq P\{L\}$, defined in Lemma 3. $C_F(H)$ and $C_F(L)$ are
matrix Cayley-Dickson algebras and $C_F(H)$ is a non-associative
subalgebra of $C_F(L)$. The algebras $C_F(H)$, $C_F(L)$ are
isomorphic as split non-associative composition algebras over the
same field $F$.

By the supposition, $H \subset L$. Then it follows from the
isomorphism of composition algebras $F\{L\}$ and $F\{H\}$ that
element $1 \neq a \in L \backslash H$ is linearly expressed
through the elements of loop $H$ in algebra $F\{H\} = C_F(H)$. Let
$FH$ be a loop algebra (without inverted commas) of loop $H$. It
follows from the definition of $''$the loop algebra$''$ $C_F(H)$
that $C_F(H) = FH/I$, where $I$ is the ideal of loop algebra
(without inverted commas) $FH$ [9]. It follows from here that in
loop algebra $FH$ element $a \in L \backslash H$ is linearly
expressed through the elements of loop $H$. Further, $H \subset
L$, therefore $FH \subseteq FL$. Then in loop algebra $FL$ element
$a \in L$ is linearly expressed through the elements of loop $H
\subset L$. But this contradicts the definition of loop algebra
$FL$, which is a free $F$-module with basis consisting of elements
of loop $L$. Consequently, the simple Moufang loop $M(F)$ has no
proper over field $F$ non-associative subloops. This completes the
proof of Lemma 6.
\smallskip\\

\textbf{Theorem 1.} \textit{Let $P$ be an algebraically closed
field. Only and only the loops $M(F)$ of the matrix Cayley-Dickson
algebra $C(F)$, where $F$ is a subfield of field $P$, are with
precise till isomorphism non-associative simple Moufang loops.
Loop $M(F)$ is quotient loop $M_0(F)/<-1>$, where $M_0(F)$
consists of all elements of  $C(F)$ with norm $1$, and the subloop
$<-1>$, generated by element $-1$, coincide with the center of
loop $M_0(F)$.}
\smallskip\\

\textbf{Proof.} If $F$ is an arbitrary subfield of $P$ then by
Lemma 4 the loop $M(F)$ is a simple non-associative Moufang loop.
Let now $Q$ be an arbitrary non-associative simple Moufang loop.
By Lemma 3 the loop $\overline{\overline Q}$ is embedded in loop
$M_0(P)$. We identify $\overline{\overline Q}$ with isomorphic
image in $M_0(P)$. Let loop $\overline{\overline Q}$ be presented
by matrices $a = \left(
\begin{array}{ll} \alpha_1 & \alpha_{12} \\ \alpha_{21} &\alpha_2
\end{array} \right)$ and let $F$ be a subfield of field $P$,
generated by all matrices elements of matrices $a$. Loop
$\overline{\overline Q}$  is a loop over field $F$. By Lemma 6
there is only one non-associative  Moufang loop over field $F$,
and namely $M_0(F)$. Therefore $\overline{\overline Q} = M_0(F)$.
By Paige's results [3] (presented at the beginning of the article)
the loop $M_0(F)$ posed only one homomorphism: $M_0(F) \rightarrow
M_0(F)/<-1> = M(F)$. Then the homomorphism $\rho:
\overline{\overline Q} \rightarrow Q$ coincides with this
homomorphism. Hence $Q = M(F)$. This completes the proof of
Theorem 1.
\smallskip\\

It is worth mentioning that in [13] a particular case of Theorem 1
is proved through other means.

 It is known that the field of complex numbers is algebraically
closed and contains as subfields all finite fields. Then from
Theorem 1 there follows the main result of article [4] about the
classification of finite non-associative simple Moufang loops,
conducted with the help of the finite groups theory.

\smallskip
TIRASPOL STATE UNIVERSITY OF MOLDOVA$$ $$

The author's home address:

Sandu Nicolae Ion

Deleanu str 1,

Apartment 60

Kishinev MD-2071, Moldova

E-mail: sandumn@yahoo.com


\begin{thebibliography}{13}
\bibitem{1} Bruck R. H. \textit{A survey of binary systems.}
Berlin-Gottingen-Heidelberg, Springer-Verlag, 1958.
\bibitem{2} Moufang R. \textit{Zur Structur von
Alternativekorpern.} Math. Ann., 1935, 110, 416 -- 430.
\bibitem{3} Paige L. J. \textit{A class of simple Moufang loops.} Proc. Amer. Math.
Soc., 1956, 7,  471 - 482. \bibitem{4} Liebeck M. W. \textit{The
classification of finite Moufang loops.} Math. Proc. Camb. Phil.
Soc., 1987, 102, 33, 33 - 47. \bibitem{5} Glauberman G. \textit{On
loops of odd order II.} J. Algebra, 1968, 8, 393 - 414.
\bibitem{6} Doro S. \textit{Simple Moufang loops.} Math. Proc. Camb. Phil.
Soc., 1978, 83, 377 - 392. \bibitem{7} Schafer R. D. \textit{An
introduction to nonassociative algebras.} New York, Academic
Press, 1966. \bibitem{8} K. A. Zhevlakov, A. M. Slin'ko, I. P.
Shestakov, A. I. Shirshov, \textit{ Rings that are nearly
associative.} Nauka, Moscow, 1978; English transl., Academic
Press, 1982.
\bibitem{9} Sandu N. I. \textit{About the embedded of Moufang
loops in alternative algebras.} To appear.
\bibitem{10} Leong F. \textit{The devil and the angel of loops.} Proc. Amer.
Math. Soc., 1976, 54, 32 - 34. \bibitem{11} Sandu N. I.
\textit{Associative nilpotent loops (Russian).} Bul. Akad.
Stiince, RSS Moldoven, 1978, 1, 28 - 33. \bibitem{12} Kleinfeld E.
\textit{On simple alternative rings.} Amer. J. Math., 1953, 75, 1,
98 - 104. \bibitem{13} Loghinov E. K. \textit{About embedding of
strongly simple Moufang loops in simple alternative algebras
(Russian).} Matem. zametki, 1993, 54, 6, 66 - 73.
\end{thebibliography}
\end{document}